\documentclass[11pt]{amsart}
\numberwithin{equation}{section}

\usepackage{amsmath,amssymb,amscd,amsxtra,latexsym,
graphicx,epsfig,euscript,amsthm,mathscinet}
\usepackage{amsmath,amsfonts,amsthm,amsopn,amssymb,enumitem}
\usepackage{graphicx}
\usepackage{float}
\usepackage{a4wide}
\usepackage{color}
\usepackage{esint}
\usepackage{palatino}
\usepackage[backend=bibtex,style=numeric,sorting=nyt,
giveninits=true,isbn=false,url=false,doi=false]{biblatex}
\renewbibmacro{in:}{\ifentrytype{article}{}{\printtext{\bibstring{in}\intitlepunct}}}
\DeclareFieldFormat*{title}{#1}
\usepackage[colorlinks]{hyperref}
\AtBeginDocument{%
  \hypersetup{%
    linkcolor=blue,%
    citecolor=green,%
  }%
}
\usepackage{cleveref}
\crefname{equation}{}{}
\crefname{figure}{}{}

\makeatletter
\@namedef{subjclassname@2020}{\textup{2020} Mathematics Subject Classification}
\makeatother

\newtheorem{theorem}{Theorem}[section]
\newtheorem{lemma}[theorem]{Lemma}

\newtheorem{proposition}[theorem]{Proposition}
\newtheorem{remark}[theorem]{Remark}

\def\e{\varepsilon}

\begin{document}

\title[Bifurcating domains for an overdetermined eigenvalue problem in cylinders]{Bifurcating domains for an overdetermined eigenvalue problem in cylinders}

\author{Yuanyuan Lian}
\address{(Y.~Lian)
  Departamento de An\'alisis matem\'atico, Universidad de Granada,
  Campus Fuentenueva, 18071 Granada, Spain}
\email{lianyuanyuan.hthk@gmail.com; yuanyuanlian@correo.ugr.es}

\author{Filomena Pacella}
\address{(F.~Pacella)
  Dipartimento di Matematica, Università di Roma Sapienza, P.le Aldo Moro 5, 00185 Roma, Italy}
\email{filomena.pacella@uniroma1.it}

\author{Pieralberto Sicbaldi}
\address{(P.~Sicbaldi)
  IMAG, Departamento de An\'alisis matem\'atico, Universidad de Granada,
  Campus Fuentenueva, 18071 Granada, Spain
  \& Aix Marseille Universit\'e - CNRS, Centrale Marseille - I2M, Marseille, France}
\email{pieralberto@ugr.es}
\thanks{{\bf Acknowledgements.} F. P. is member of GNAMPA-INdAM and has been supported by PRIN 2022AKNS E4 - Next Generation EU.
Y. L. and P. S. have been supported by the Grants PID2020-117868GB-I00 and PID2023-150727NB-I00 of the MICIN/AEI.
P. S. has been supported also by the \emph{IMAG-Maria de Maeztu} Excellence Grant CEX2020-001105-M funded by the MICIN/AEI.}



\maketitle

\noindent

\noindent
\begin{abstract}
We study an overdetermined eigenvalue problem for domains $\Omega$ contained in the half-cylinder $\Sigma=\omega \times (0, +\infty)$, based on a bounded regular domain $\omega \subset \mathbb{R}^{N-1}$. It is easy to see that in any bounded cylinder $\Omega_{t}=\omega \times (0, t)$, $t > 0$, the eigenvalue problem admits a one-dimensional positive eigenfunction which satisfies the overdetermined boundary conditions. The aim of the paper is to construct other domains $\Omega\subset \Sigma$ for which there exists a positive eigenfunction that is a solution of the overdetermined problem. This is achieved by showing that branches of such domains bifurcate from the ``trivial'' domains $\Omega_{t_j}$ at the values $t_{j} = \frac{\pi}{2\sqrt{\sigma_j}}$ where $\sigma_j$ ($j\geq 1$) is a simple Neumann eigenvalue of the Laplace operator on $\omega \subset \mathbb{R}^{N-1}$. The solutions can be reflected with respect to $\omega$ to generate nontrivial solutions in a cylinder.

\end{abstract}

\section{Introduction}
\label{Section 1}

Let $\Omega$ be a domain of $\mathbb{R}^N$ ($N\geq 2$). The classical overdetermined eigenvalue problem in $\Omega$ is given by:
\begin{equation}\label{eigenvalue}
	\begin{cases}
		\Delta {u}+\lambda u=0 &\mbox{in $\Omega$},\\
		u=0 &\mbox{on $\partial \Omega$},\\
		\partial_\nu u = c &\mbox{on $\partial \Omega$},
	\end{cases}
\end{equation}
where $\nu$ is the exterior unit normal on $\partial \Omega$ and $c$ is a constant. Here we suppose that $\lambda$ is the first eigenvalue of the Laplacian in $\Omega$, so the solution $u$ will be positive. The constant $\lambda$ is given by
\begin{equation*}
  \lambda_1(\Omega)=\min_{v\in H_0^1(\Omega)} \frac{\int_{\Omega} |\nabla v|^2}{\int_{\Omega}v^2}.
\end{equation*}
If $\Omega$ is a bounded and regular domain, according to the celebrated work by Serrin \cite{MR333220}, the domain $\Omega$ must be a ball and the solution $u$ is radial symmetric. For unbounded domains, Berestycki, Caffarelli and Nirenberg \cite{MR1470317} considered
\begin{equation}\label{e.Laplacian}
  \begin{cases}
  \Delta u+f(u) =0& \mbox{in $\Omega$, }\\
  u>0 &\mbox{in $\Omega$, }\\
  u=0 &\mbox{on $\partial\Omega$, }\\
  \partial_{\nu} u=c &\mbox{on $\partial\Omega$, }
  \end{cases}
\end{equation}
where $f$ is a locally Lipschitz function and conjectured that if $\Omega$ is a smooth domain with $\mathbb{R}^N \setminus \bar{\Omega}$ connected, the existence of a solution $u$ to the previous problem implies that $\Omega$ is either a ball, a half-space, a generalized cylinder ($B^k \times \mathbb{R}^{N-k}$ where $B^k$ is a ball in $\mathbb{R}^k$) or the complement of one of them.
This conjecture has been disproved by the third author \cite{MR2592974} in dimensions $N\geq 3$, using a bifurcation argument. Precisely, he showed that there exist periodic domains given by perturbations of the straight cylinder such that problem \eqref{eigenvalue} admits a positive solution. Since then many counterexamples had been built by bifurcation arguments, see \cite{MR2854185, MR3417183, MR4046014, MR4484836, MR4649188}.

\medskip

If we change the ambient space $\mathbb{R}^N$ by a certain open set $\Sigma \subset \mathbb{R}^N$ and we suppose that the domain $\Omega$ is contained in $\Sigma$, then we have a {\it relative} overdetermined eigenvalue problem in $\Sigma$, i.e.,
\begin{equation}\label{main_1}
	\begin{cases}
		\Delta {u}+\lambda u=0 &\mbox{in $\Omega$},\\
		u=0 &\mbox{on $\Gamma_\Omega$},\\
        \partial_{\eta} u = 0 &\mbox{on $\widetilde \Gamma_\Omega$},\\
        \partial_{\nu}u =c  &\mbox{on $\Gamma_\Omega$},
	\end{cases}
\end{equation}
where $\Gamma_\Omega = \partial \Omega \cap \Sigma$, $\widetilde \Gamma_\Omega = \partial \Omega \cap \partial \Sigma$, $\nu$ is the exterior unit normal on $\Gamma_\Omega$ and $\eta$ is the exterior unit normal on $\widetilde \Gamma_\Omega$. Observe that the overdetermined condition appears only on $\Gamma_\Omega$, that is usually called the {\it relative} or {\it free} boundary of $\Omega$. We consider again that $\lambda$ is the first eigenvalue of the Laplacian in $\Omega$ and then the corresponding solution $u$ is supposed to be positive. A natural question is to determine the domains $\Omega \subset \Sigma$ for which \eqref{main_1} admits a solution. It is worth to point out that, by \cite{MR3417691}, if the intersection between $\Gamma_{\Omega}$ and $\widetilde \Gamma_{\Omega}$ is transversal, then $\Gamma_{\Omega}$ and $\widetilde \Gamma_{\Omega}$ must be orthogonal.

\medskip

In the relative setting, the overdetermined torsion problem, i.e,
\begin{equation}\label{torsion}
	\begin{cases}
		\Delta {u}+1=0 &\mbox{in $\Omega$},\\
		u=0 &\mbox{on $\Gamma_\Omega$},\\
        \partial_{\eta} u = 0 &\mbox{on $\widetilde \Gamma_\Omega$},\\
        \partial_{\nu}u =c  &\mbox{on $\Gamma_\Omega$},
	\end{cases}
\end{equation}
has been studied in \cite{MR4109828, MR4185827, MR4668558, MR4689375, MR4887781}. In \cite{MR4109828, MR4185827}, the second author and Tralli considered \eqref{torsion} when $\Sigma$ is a cone. If $\Sigma$ is a cylinder, \eqref{torsion} has been considered in \cite{MR4668558, MR4689375, MR4887781}. In \cite{MR4668558}, the authors introduced the definition of relative Cheeger set (\cite[Definition 3.1]{MR4668558}) for a domain $\Omega \subset \Sigma$ when $\Sigma$ is a general Lipschitz unbounded domain. They showed that if $\Sigma$ is convex, any bounded domain for which \eqref{torsion} admits a solution coincides with its relative Cheeger set. In \cite{MR4689375, DeMMP2025}, a more general relative overdetermined problem has been studied and the stability analysis suggests that nontrivial solutions to problem \eqref{torsion} could exist. Later, the second author, Ruiz and the third author \cite{MR4887781} proved that there are nontrivial domains for which \eqref{torsion} admits a solution.

\medskip

In this paper, we consider problem \eqref{main_1} in
\begin{equation}
   \Sigma = \Sigma_\omega := \{(x', x_N) \in \mathbb{R}^{N} \ | \ x' \in \omega, \ x_N \in (0, + \infty)\},
\end{equation}
where $\omega\subset \mathbb{R}^{N-1}$ is a bounded regular domain. Note that in the domains $\Omega_t=\omega\times (0,t)$, $\forall t>0$, problem \eqref{main_1} always admits $1$-dimensional positive eigenfunctions, i.e., eigenfunctions which depend only on the $x_N$ variable. We refer to $\Omega_t$ as a trivial domain. Our aim is to build new nontrivial domains $\Omega \subset \Sigma$ for which \eqref{main_1} admits a positive solution. Such domains $\Omega$ will be the hypographs of $C^2$ functions $v : \omega \to \mathbb{R}^+$, i.e.,
\begin{equation*}
    \Omega = \Omega_v := \{(x', x_N) \in \mathbb{R}^N \ | \ x' \in \omega, \ x_N \in (0, v(x'))\},
\end{equation*}
where $\Gamma_{\Omega_v}$ is just the cartesian graph of $v$ over $\omega$. In this case we write for short:
\begin{equation*}
  \Gamma_{v}:=\Gamma_{\Omega_v},\quad \widetilde \Gamma_{v}:=\widetilde \Gamma_{\Omega_v},
\end{equation*}
where $\widetilde \Gamma_{\Omega_v}=\partial \Omega_v \cap \partial \Sigma_{\omega}$. Then we consider positive solutions of the following overdetermined eigenvalue problem:
\begin{equation}\label{main}
	\begin{cases}
		\Delta {u}+\lambda u=0 &\mbox{in $\Omega_v$},\\
        u=0 &\mbox{on $\Gamma_v$},\\
        \partial_\eta u = 0 &\mbox{on $\widetilde \Gamma_v$},\\
		\partial_{\nu}u =c  &\mbox{on $\Gamma_v$}.
	\end{cases}
\end{equation}
In view of the orthogonal condition between $\Gamma_{\Omega}$ and $\widetilde \Gamma_{\Omega}$ recalled before, we point out that on $\partial \Gamma_v$ the two normal vectors are orthogonal if $\partial_{\eta} v =0$ on $\partial \omega$ \cite{MR4887781}.

\medskip

To state our result, we consider the sequence of eigenvalues of the Laplacian operator on $\omega$ with Neumann boundary condition:
\begin{equation*}
  0 = \sigma_0 < \sigma_1  \leq \sigma_2 \cdots.
\end{equation*}

Our result is the following:

\begin{theorem} \label{Tmain}
Assume that $\sigma_j$ is simple for some $j \geq 1$. Then there exists a smooth curve:
\begin{equation*}
  (-\varepsilon, \varepsilon) \ni s \mapsto v(s) \in C^{2,\alpha}(\bar\omega)
\end{equation*}
such that:
\begin{enumerate}
	\item $v(0)$ is the constant function $\frac{\pi}{2\sqrt{\sigma_j}}$ ;
	\item $v(s)$ is a nonconstant function for $s \neq 0$;
	\item For every $s\in (-\varepsilon, \varepsilon)$, the problem \eqref{main} admits a positive eigenfunction $u=u_{v(s)}$ with eigenvalue $\lambda=\lambda_{v(s)}$ in the domain $\Omega=\Omega_{v(s)}$, where the constant $c$ depends on $s$.
    \item $\Gamma_{v(s)}$ intersects $\partial \Sigma_\omega$ orthogonally.
\end{enumerate}
\end{theorem}
We remark here that, for most bounded $C^2$ regular domains $\omega \subset \mathbb{R}^{N-1}$, all Neumann eigenvalues are simple \cite[Chapter 6, Example 6.4]{MR2160744}. In fact, our result does not depend on the particular cylinder. The important thing is that at least one Neumann eigenvalue $\sigma_j$ is simple. The shape of $\omega$ is irrelevant.

The proof of \Cref{Tmain} uses bifurcation analysis. In the context of overdetermined problems, it was first used in \cite{MR2592974}. In the case of the relative overdetermined torsion problem in cylinder it has been used in \cite{MR4887781}. With respect to this last work, one of the main differences is that in our work the linearized problem is degenerate, i.e. has a kernel (spanned by the first eigenfunction) and this fact must be treated with care. Moreover the computations to use the bifurcation results are more involved. The two existence results, i.e. that of this paper and that of \cite{MR4887781} can probably be generalized to other overdetermined problems with general equation of the form $\Delta u + f(u)=0$, and according to the degeneracy or not of its linearization the strategy to obtain the result should be that of this paper or that of \cite{MR4887781}. In this sense, the two problems are the prototypes for bifurcations of overdetemined problems in cylinders.

\medskip

We finally remark that any solution of a general semilinear relative Dirichlet problem in such $\Omega_{v}\subset \Sigma_{\omega}$:
\begin{equation}\label{eq:1.7}
\begin{cases}
\Delta u+f(u)=0 &~~\mbox{in }~~\Omega_{v}, \\
u>0 &~~\mbox{in }~~\Omega_{v}, \\
u=0 &~~\mbox{on }~~\Gamma_{v}, \\
\partial_{\eta}u=0 &~~\mbox{on }~~\widetilde{\Gamma}_{v},
\end{cases}
\end{equation}
where $f$ is Lipschitz continuous, can be reflected with respect to the base $\omega$ of $\Sigma_{\omega}$, providing a solution of the analogous problem in the domain:
\begin{equation*}
  D_{v}=\{(x',x_{N})\in \mathbb{R}^{N} \mid x'\in\omega, -v(x')<x_N<v(x') \}
\end{equation*}
whose relative boundary is the union of $\Gamma_{v}$ and $\Gamma_{-v}$. Therefore, our result also gives nontrivial solutions for the relative overdetermined eigenvalue problem for domains $D_{v}$ which are contained in the whole cylinder $Z_{\omega}=\omega\times(-\infty,+\infty)$. Since any domain $D_{v} \subset Z_{\omega}$ is symmetric with respect to $\omega$ and convex in the $x_N$-direction, by applying the moving plane method, it is easy to see that any solution of:
\begin{equation}\label{eq:1.8}
\begin{cases}
\Delta u+f(u)=0 &~~\mbox{in }~~D_{v}, \\
u>0 &~~\mbox{in }~~D_{v}, \\
u=0 &~~\mbox{on }~~\Gamma_{v}\cup \Gamma_{-v}, \\
\partial_{\eta} u=0 &~~\mbox{on }~~\partial D_{v}\cap \Sigma_{\omega},
\end{cases}
\end{equation}
is even in the $x_N$-variable and so gives a solution of \eqref{eq:1.7}. Thus, to study the overdetermined problem in domains $\Omega_{v} \subset \Sigma_{\omega}$ or domains $D_{v} \subset Z_{\omega}$ is equivalent.

The paper is organized as follows. In \Cref{sec:rel_dic_prob}, a normal derivative operator $F$ associated to each hypograph is defined and studied. In \Cref{sec:linearization} and \Cref{sec:linearized_operator_H_t}, we calculated the Fr\'{e}chet derivative of the operator $F$ at any constant function and analyze its properties. In \Cref{sec:proof_Tmain}, we prove \Cref{Tmain} by using the Crandall-Rabinowitz Bifurcation Theorem.

\section{The mixed boundary eigenvalue problems and the normal derivative operator}
\label{sec:rel_dic_prob}
In this section, we first recall some basic notions about the eigenvalue problem with mixed boundary conditions and then introduce the $1$-dimensional associated eigenfunction. Next, we will show that the first eigenvalue $\lambda_v$ in the domain $\Omega_v$, with its eigenfunction $u_v$, depend smoothly on $v$. In the following, we use $x=(x', x_N)=(x_1,...,x_N)$ to denote the coordinates in $\mathbb{R}^N$. In addition, $u_i$ denotes $\partial u/\partial x_i$ and $u_{ij}$ denotes $\partial^2u/\partial x_i\partial x_j$ ($1\leq i,j\leq N$) etc. We also use the common Einstein notation dropping the symbol of sum, that is understood when indices are repeated.

\medskip

Fix $\alpha \in (0,1)$.  For $k \in \mathbb{N}$, define:
\begin{equation} \label{space}
X_k= \{ v \in C^{k,\alpha}(\bar \omega):\  \partial_{\eta} v=0 \mbox{ on }\partial \omega\}
\end{equation}
and
\begin{equation} \label{space2}
X_k^+=\left\{v\in X_k: \inf_{x'\in \omega} v(x')>0\right\}, \quad \widetilde{X}_k= \{ v \in X_k: \ \int_{\omega} v=0\}
\end{equation}
endowed with the usual $C^{k, \alpha}$ norm which is denoted as $\| \cdot \|_{C^{k,\alpha}}$.
Note that any function $v \in X_k$ can be decomposed as $v = \tilde{v} + h$, where $\tilde{v} \in \widetilde{X}_k$ and $h \in \mathbb{R}$. Indeed,
\begin{equation*}
  h = \frac{1}{|\omega|} \int_{\omega} v, \quad \tilde{v}=v-h.
\end{equation*}

We have the following existence result concerning the eigenvalues and eigenfunctions of the Laplace operator with mixed boundary conditions.

\medskip

For any $v \in X_2^+$, consider the eigenvalue problem:
\begin{equation}\label{dir}
	\begin{cases}
		\Delta {u}+\lambda u=0 &\mbox{in $\Omega_{v},$}\\
		u=0 &\mbox{on $\Gamma_{v},$}\\
		\partial_\eta u = 0 &\mbox{on $\widetilde \Gamma_{v}.$}
	\end{cases}
\end{equation}
It admits a sequence of eigenvalues $\{\lambda_j\}_{j\in \mathbb{N}^+}\subset \mathbb{R}^+$ with
\begin{equation*}
0<\lambda_1<\lambda_2\leq \lambda_3\leq \cdots\leq \lambda_j\to +\infty
\end{equation*}
and  a sequence of corresponding $L^2$-normalized eigenfunctions $\{u_j\}_{j\in \mathbb{N}^+}\subset H_0^1(\Omega_v \cup \widetilde \Gamma_v)$ such that $u_1>0$ and $\lambda_1$ is simple in $\Omega_v$ and $\{u_j\}$ is an orthonormal basis of the space $L^2(\Omega_v)$. Here, $H_0^1(\Omega_v \cup \widetilde \Gamma_v)$ is the subspace of $H^1(\Omega_v)$ made by functions with trace vanishing on $\Gamma_v$.

\medskip

The proof of the existence of eigenvalues with mixed boundary condition can be obtained by the standard spectral theory for symmetric bilinear forms (as in \cite[Theorem 6.74 and Theorem 6.76]{MR3497072}). The key is that the Poincar\'{e} inequality holds because of $u=0$ on $\Gamma_v$. This is similar to the problem with Dirichlet boundary condition \cite[Theorem 8.8]{MR3497072} (see also \cite{MR1814364}). In \cite{MR4319030}, Damascelli and the second author considered the existence of a more general eigenvalue problem with mixed boundary conditions.

\medskip

If $v(x')\equiv t$, $\forall x' \in \omega$ and $t\in \mathbb{R}^+$, then $\Omega_v=\omega\times (0,t)$ and we use the notations $\Omega_t$, $\Gamma_t$ and $\widetilde \Gamma_t$ for the corresponding $\Omega_v$, $\Gamma_v$ and $\widetilde \Gamma_v$. Then the problem \eqref{dir} admits the following $1$-dimensional first eigenfunction:
\begin{equation}\label{so:U_t}
   u_t(x', x_N):=\cos \left(\frac{\pi}{2t}x_N\right)
\end{equation}
corresponding to the first eigenvalue:
\begin{equation*}
\lambda_t=\left(\frac{\pi}{2t}\right)^2.
\end{equation*}
In fact, we can obtain $u_t$ and $\lambda_t$ from the following ODE with both Dirichlet and Neumann boundary conditions:
\begin{equation}\label{eq:U_t}
	\begin{cases}
		u''+\lambda u=0 &~~\mbox{ in }~~(0,t),\\
        u>0 &~~\mbox{ in }~~(0,t),\\
		u(t)=0,\\
		u'(0) =0.
	\end{cases}
\end{equation}

In this paper, we denote $u_t(x', x_N)$ by $u_t(x_N)$ for simplicity and often refer to $\Omega_t$ as the trivial domains. Note that $u_t>0$ in $\Omega_t$. Hence, $\lambda_t$ is simple and $u_t$ is the unique eigenfunction (up to linear dependence).

In addition, by a direct calculation, $u_t$ is also a solution of the overdetermined problem \eqref{main} with $\lambda=\lambda_t$ and $c=-\pi/2t$. In this paper, we will carry out the bifurcation from $u_t$ and $\lambda_t$ for suitable values of $t$.

\medskip

In the following, we show that for the domain $\Omega_v$, the first eigenvalue $\lambda_v$ and its eigenfunction $u_v$ of \eqref{dir} depend on $v$ in a smooth way. The relation between $u_v$ and $v$ is not explicit since $u_v$ is determined by the domain $\Omega_v$ which is influenced by $v$. Thus, we rephrase problem \eqref{dir} in order to understand the relation between the solution $u_v$ and the domain $\Omega_v$ in a more explicit sense.

Let us define the following diffeomorphism $Y_v: \Omega_1 \to \Omega_{v}$,
\begin{equation} \label{Y}
	x=(x', x_N)=Y_v(y)=Y_v (y', y_N) = \left(y', v (y') y_N \right),
\end{equation}
which is valid for any $v \in X_k^+$. The coordinates we consider from now on are $(y',y_N) \in \Omega_1$ and in this coordinate system, the metric $g$ is given by the pullback of the Euclidean metric by $Y_v$ which can be written as
\begin{equation} \label{metrica}
	g_{ij}=
\left(
\begin{matrix}
\delta^{ij}+v_iv_jy_{N}^2 & vv_iy_N\\
vv_iy_N & v^2
\end{matrix}
\right),\quad i,j=1,\cdots,N-1.
\end{equation}
Thus, the determinant of $g$ is
\begin{equation*}
  \det g=v^2.
\end{equation*}
The inverse of $g$ is given as
\begin{equation}\label{g-inverse}
  g^{ij}=
\left(
\begin{matrix}
\delta^{ij} & -\frac{v_iy_N}{v}\\
-\frac{v_iy_N}{v} & \frac{1+|Dv|^2y_N^2}{v^2}
\end{matrix}
\right).
\end{equation}

Let
\begin{equation*}
  \phi= Y_v^* u_v \quad \left(\mbox{i.e.,}~~\phi(y)=u_v\left(Y_v(y)\right)=u_v(x) \right).
\end{equation*}
Then, $u_v$ is a solution of problem \eqref{dir} if and only if $\phi(y',y_N)$ solves:
\begin{equation}\label{formula-new}
	\begin{cases}
     \Delta_g \phi + \lambda \phi  = 0 & \mbox{in $\Omega_1$},\\
	 \phi = 0 &\mbox{on $\Gamma_1$},\\
	 \partial_{ \eta}  \phi =  0 &\mbox{on $\tilde \Gamma_1$},
    \end{cases}
\end{equation}
where $\Delta_{g}$ is the Laplace-Beltrami operator with respect to the metric $g$. Denote
\begin{equation*}
  h(y') = \frac{1}{v(y')},
\end{equation*}
then we have
\begin{eqnarray*}
		\Delta_g \phi & = & \sum_{i=1}^{N-1} \left[ \phi_{ii} + \frac{y_N}{h}\, \left( 2 \phi_{iN}\, h_i  + \phi_{N}\, h_{ii} + \phi_{NN}\, h_i^2\, \frac{y_N}{h}\right) \right] + \phi_{NN}\, h^2.
\end{eqnarray*}

\medskip

Next, inspired by \cite[Proposition 3.2]{MR2521426}, we show that the first eigenvalue $\lambda_v$ and its eigenfunction $u_v$ depend on $v$ in a smooth way. Define
\begin{equation*}
  Z = \{ \phi \in C^{2, \alpha}(\bar \Omega_1): \ \phi=0 \mbox{ on } \Gamma_1, \ \partial_\eta \phi=0 \mbox{ on } \widetilde{\Gamma}_1 \}.
\end{equation*}

\begin{proposition}
\label{dirichlet}
For any $v \in X_2^+$, there exist a unique positive solution $\phi_v \in Z$ of the problem \eqref{formula-new} with $\|\phi_v\|_{L^2(\Omega_1)}=\sqrt{|\omega|/2}$ and a unique corresponding eigenvalue $\lambda_v$ (equivalently, there exist a unique positive solution $u_v$ of \eqref{dir} with $\lambda_v$ and $\|u_v/\sqrt{v}\|_{L^2(\Omega_v)}=\sqrt{|\omega|/2}$). In addition the map:
\begin{equation} \label{defPhi}
\mathcal{T}: X_2^+ \to Z\times \mathbb{R}^+, \quad  \mathcal{T}(v):= (\phi_v, \lambda_v)
\end{equation}
is $C^{\infty}$.
\end{proposition}
\begin{remark}\label{re2.3}
The reason why we require $\|\phi_v\|_{L^2(\Omega_1)}=\sqrt{|\omega|/2}$
where $|\omega|$ is the $(N-1)$-dimensional measure of $\omega$
is due to the following fact:
\begin{equation*}
\|\Phi_t\|_{L^2(\Omega_1)}=\sqrt{|\omega|/2},
\end{equation*}
where $\Phi_t(y):=u_t(x)$ and $x=Y_t(y)$.
\end{remark}

\begin{proof}
Given $v_0\in X_2^+$, as mentioned previously, there exists a unique $(u_0,\lambda_0)$ solution of \eqref{dir} with $\|u_0 / \sqrt{v} \|_{L^2(\Omega_v)}=\sqrt{|\omega|/2}$. Equivalently, there exist a unique $(\phi_0,\lambda_0)=(\phi_{v_0},\lambda_0)\in Z\times \mathbb{R}^+$ that solves \eqref{formula-new} with $\|\phi_0\|_{L^2(\Omega_1)}=\sqrt{|\omega|/2}$, where $\phi_0(y) = u_0(x)$ and $x=Y_{v_0}(y)$. In the following we just show the smooth dependence of $(\phi,\lambda):=(\phi_v,\lambda_v)$ on $v$.

Let
\begin{equation*}
Z_{\bot}=\left\{\phi\in Z:  \int_{\Omega_1} \phi \cdot\phi_0=0 \right\},\quad C^{\alpha}_{\bot}=\left\{\phi\in C^{\alpha}(\overline \Omega_1): ~
\int_{\Omega_1}\phi \cdot\phi_0=0\right\}.
\end{equation*}
Define the mapping $T: X_2^+\times Z_{\bot}\to C^{\alpha}_{\bot}$ as
\begin{equation*}
T(v,\phi)=L_{v_0} \phi +\lambda_0 \phi+\left(L_{v}-L_{v_0}+\mu\right)\left(\phi_0+\phi\right),
\end{equation*}
where $L_v:=\Delta_g $ and $g$ is the induced metric on $\Omega_1$ by $Y_v$ (see \eqref{metrica}). The constant
$\mu$ is given by
\begin{equation*}
\mu=-\int_{\Omega_1} \phi_0(L_{v}-L_{v_0})(\phi_0+\phi),
\end{equation*}
so that $T(v,\phi)\in C^{\alpha}_{\bot}$. Clearly, $T$ is a $C^{\infty}$ mapping.

In the following, we try to use the implicit function theorem to prove the smooth dependence. First, we note that \begin{equation}\label{e4.2}
T(v_0, \phi_0) =0.
\end{equation}
Then by a direct calculation,
\begin{equation*}
L(\psi):=D_{\phi}T_{(v_0,\phi_0)}(\psi)= L_{v_0} \psi + \lambda_0 \psi.
\end{equation*}

Consider the following problem:
\begin{equation}\label{e2.4}
	\begin{cases}
		L_{v_0} \psi + \lambda_0 \psi= f &\mbox{in $\Omega_{1},$}\\
		\psi =0 &\mbox{on $\Gamma_{1},$}\\
		\partial_\eta \psi = 0 &\mbox{on $\tilde \Gamma_{1}.$}
	\end{cases}
\end{equation}
By the classical spectral theory (e.g., \cite[Theorem 6.66, Page 402]{MR3497072}), \eqref{e2.4} has a nonzero solution if and only if
\begin{equation*}
\int_{\Omega_1}f\cdot \phi_0=0,
\end{equation*}
which is exactly the definition of $C^{\alpha}_{\bot}$. Hence, $L$ is an invertible operator from $Z_{\bot}$ to $C^{\alpha}_{\bot}$. By the implicit function theorem, for any $v\in \mathcal{N}$, where $\mathcal{N}$ is a neighbourhood of $v_0$, there exists a unique $\phi\in Z_{\bot}$ such that
\begin{equation}\label{e4.3}
T(v,\phi)=0,
\end{equation}
and the dependence of $\phi$ on $v$ is smooth. By its definition, $\mu$ depends smoothly on $v$ as well. The \eqref{e4.3} is equivalent to
\begin{equation*}
L_v(\phi_0+\phi)+(\lambda_0+\mu)(\phi_0+\phi)=0.
\end{equation*}
Thus, $\phi_0+\phi$ is a solution of \eqref{formula-new} with the eigenvalue $\lambda_v:=\lambda_0+\mu$. Finally, define
\begin{equation*}
\phi_{v}=c(\phi_0+\phi),
\end{equation*}
where $c$ is a constant such that $\|\phi_v\|_{L^2(\Omega_1)}=\sqrt{|\omega|/2}$.
\end{proof}

\medskip

For a fixed $t_*>0$ (which will be chosen later), let us consider the interval $I_{\delta}=(t_*-\delta,t_*+\delta)$ for some fixed small $\delta>0$ such that $t_*-\delta>0$. Then we consider the ball $B_{\rho}$ in $\widetilde{X}_2$ with center in $0$ and radius $\rho>0$ such that
\begin{equation*}
  t+w(x')>0\quad \forall t\in I_{\delta}, ~~ \forall x'\in \omega ~~\mbox{ and }~~\forall w \in B_{\rho}.
\end{equation*}

The previous proposition allows us to define the operator:
\begin{equation*}
  F: I_{\delta}\times B_{\rho} \to \widetilde{X}_1
\end{equation*}
by
\begin{equation}\label{F}
F (t,w) (x') =  \displaystyle  \frac{\partial u_{t+w}}{\partial {\nu}} (x', t+w(x')) - \frac{1}{|\omega|} \, \int_{\omega} \, \frac{\partial u_{t+w}}{\partial \nu},
\end{equation}
where $u_{t+w}$ is the solution of \eqref{dir} with $v=t+w$. By \Cref{dirichlet}, $F$ is a smooth operator. Observe that $F$ is well defined since on $\partial \Gamma_{t+w}$,
\begin{equation*}
  \frac{\partial}{\partial \eta} \left (\frac{\partial u_{t+w}}{\partial \nu} \right )=\frac{\partial}{\partial \nu} \left (\frac{\partial u_{t+w}}{\partial \eta} \right )=0.
\end{equation*}
Clearly the zeros of $F$ correspond to the solutions of \eqref{main} in $\Omega_{t+w}$. We plan to find nontrivial zeroes of $F$ as a local bifurcation of a family of trivial solutions. Indeed, if $w {\equiv} 0 $, the unique solution of \eqref{dir} is given by \eqref{so:U_t}. Clearly the Neumann derivative is constant so that $F(t,0) =0 $ for any constant $t>0$.

\section{The linearization of the normal derivative operator}
\label{sec:linearization}

In this section, we will calculate the Fr\'{e}chet derivative of $F$ with respect to $w$ at $w=0$ for any $t\in I_{\delta}$. We denote such derivative by $H_t$, i.e.,
\begin{equation}\label{Ht}
  H_t=D_w F(t,w)|_{w=0}.
\end{equation}

\medskip

Before we calculate the Fr\'{e}chet derivative $H_t$, we introduce the following lemma:
\begin{lemma}\label{le3.1}
For any $t>0$ and any $w\in \widetilde{X}_2$, there exists a unique solution $\hat \psi$ (orthogonal to $u_t$ in \eqref{so:U_t}) to the problem
\begin{equation}\label{psi_h2}
  \begin{cases}
  \Delta \hat{\psi}+\lambda_t \hat{\psi}=\, 0  ~~&\mbox{ in }~~\Omega_t,\\
      \hat{\psi}=\, \frac{\pi}{2t^2} w~~&\mbox{ on }~~\Gamma_t,\\
     \partial_{\eta} \hat{\psi}=\, 0~~&\mbox{ on }~~\widetilde \Gamma_t.
  \end{cases}
\end{equation}
\end{lemma}

\begin{remark}\label{re3.1}
We remark here that we have identified $\Gamma_{t}$ with $\omega$ in order to write the second equality.
\end{remark}

\begin{proof}
Let
\begin{equation*}
\psi =\hat{\psi}- \frac{\pi}{2t^2} w.
\end{equation*}
Then $\hat{\psi}$ is the solution of \eqref{psi_h2} if and only if $\psi$ is the solution of
\begin{equation}\label{e3.1}
  \begin{cases}
  \Delta \psi+\lambda_t \psi=\, f:= - \frac{\pi}{2t^2}\left(\Delta w+\lambda_t w\right)  ~~&\mbox{ in }~~\Omega_t,\\
      \psi=\, 0~~&\mbox{ on }~~\Gamma_t,\\
     \partial_{\eta} \psi=\, 0~~&\mbox{ on }~~\widetilde \Gamma_t.\\
  \end{cases}
\end{equation}
The third equality holds since $w\in \widetilde{X}_2$.

According to the classical spectral theory (e.g. \cite[Theorem 6.66, Page 402]{MR3497072}), \eqref{e3.1} has a nonzero solution if and only if $\langle f,u_t\rangle=0$, i.e.,
\begin{equation*}
\int_{\Omega_t}\left(\Delta w+\lambda_t w\right)u_t=0.
\end{equation*}
By integration by parts and considering the definitions of $w$ and $u_t$, we have
\begin{equation*}
  \begin{aligned}
\int_{\Omega_t}\left(\Delta w+\lambda_t w\right)u_t=&
-\int_{\Gamma_t} w u_t'+\int_{\Omega_t}\left(\Delta u_t+\lambda_t u_t\right)w\\
    =&-\int_{\Gamma_t} w u_t'=-u_t'(t)\int_{\Gamma_t} w=0,
  \end{aligned}
\end{equation*}
where we have used the fact that the integral average of $w$ over $\Gamma_t$ (i.e., $\omega$) is $0$.

Therefore, \eqref{e3.1} has a nonzero solution. By choosing the solution $\psi$ such that $\langle \psi, u_t\rangle=0$, we conclude that this $\psi$ is unique. Then $\hat{\psi}=\psi +\frac{\pi}{2t^2} w$ is the solution of \eqref{psi_h2}. Note that
\begin{equation*}
\langle w, u_t\rangle=\int_{\Omega_t}w u_t dx=\left(\int_{\omega}w(x')dx'\right)\cdot \left(\int_{0}^{t}u_t(x_N)dx_N\right)=0.
\end{equation*}
Hence, $\hat{\psi}$ is orthogonal to $u_t$.
\end{proof}

Next, we will show a property of $\hat \psi$. Note that the normal derivative of $\hat \psi$ on $\Gamma_{t}$ (i.e., $\partial_\nu \hat \psi$) is $\hat \psi_N (x',t)$ (i.e., $\partial \psi (x',t)/\partial x_N$) which is a function of the variable $x' \in \omega$ with $C^{1,\alpha}$ regularity.

\begin{lemma}\label{le.new}
Let $\hat \psi$ be the solution of \eqref{psi_h2}. Then,
\begin{equation}\label{e.new3}
  \int_{\omega} \partial_{\nu} \hat{\psi}=\int_{\omega} \hat{\psi}_N=0.
\end{equation}
\end{lemma}

\begin{proof}
Let $\{\xi_k\}_{k\geq 0}$ be the sequence of normalized eigenfunctions (i.e., $\|\xi_k\|_{L^2(\omega)}=1$) of the Laplacian operator on $\omega \subset \mathbb{R}^{N-1}$ with Neumann boundary conditions and $\sigma_k$ ($k\geq0$) are the corresponding sequence of eigenvalues. That is,
\begin{equation*}
  \begin{cases}
  \Delta \xi_k+\sigma_k \xi_k=0\quad &~~\mbox{ in }~~\omega,\\
  \partial_{\eta} \xi_k=0\quad &~~\mbox{ on }~~\partial \omega.
  \end{cases}
\end{equation*}
Obviously, $\sigma_0=0$ and any corresponding eigenfunction is constant. By multiplying $\xi_k$ on both sides of above equation and taking the integral by parts, we derive that
\begin{equation}\label{e.new2}
\int_{\omega} |\nabla \xi_k|^2=\sigma_k \int_{\omega} \xi_k^2=\sigma_k.
\end{equation}
\medskip
Considering the Sobolev space
\begin{equation*}
  \widetilde H^1(\omega)=\{w\in H^1(\omega),\quad \int_{\omega} w=0\},
\end{equation*}
we have that any function $w\in \widetilde H^1(\omega)$ can be written as
\begin{equation}\label{e.new}
  w(x')=\sum^{\infty}_{k=1} a_k \xi_k(x'),
\end{equation}
where $a_k\in \mathbb{R}$ are the Fourier coefficients. In particular, any $w\in \widetilde{X}_2$ can be written as above. In addition, $\hat{\psi}$ can be written as
\begin{equation*}
  \hat{\psi}(x',x_N)=\sum^{\infty}_{k=0} a_k g_k(x_N) \xi_k(x')
  =a_0g_0+\sum^{\infty}_{k=1} a_k g_k(x_N) \xi_k(x')
\end{equation*}
for some functions $g_k$. Inserting the above formula into \eqref{psi_h2}, for $k\geq 1$, we have
\begin{equation}\label{psi_h3}
  \begin{cases}
     g_k''-\sigma_kg_k+\lambda_t g_k= 0  ~~\mbox{ in }~~[0,t],\\
     g_k(t)= \frac{\pi}{2t^2},\\
     g_k'(0)= 0,\\
  \end{cases}
\end{equation}

We claim that $a_0=0$. Reasoning by contradiction, assume that $a_0\neq 0$. Observe that $\hat{\psi}=\frac{\pi}{2t^2} w$ on $\Gamma_t$ (see \eqref{psi_h2}), i.e.,
\begin{equation*}
  \hat{\psi} (x',t)=\frac{\pi}{2t^2}\cdot w(x')
  =\frac{\pi}{2t^2}\cdot \sum^{\infty}_{k=1} a_k \xi_k(x')
  =a_0g_0(t)+\sum^{\infty}_{k=1} a_k g_k(t) \xi_k(x').
\end{equation*}
Then we have $a_0g_0(t)=0$ and thus $g_0(t)=0$. So for $k=0$, \eqref{psi_h3} becomes:
\begin{equation}\label{psi_h4}
  \begin{cases}
     g_0''+\lambda_t g_0= 0  ~~\mbox{ in }~~[0,t],\\
     g_0(t)= 0,\\
     g_0'(0)= 0.\\
  \end{cases}
\end{equation}
It means that
\begin{equation*}
  g_0=u_t.
\end{equation*}
Then we have
\begin{equation*}
  \hat{\psi} (x',x_N)=a_0u_t+\sum^{\infty}_{k=1} a_k g_k(x_N) \xi_k(x'),
\end{equation*}
which contradicts that $\hat{\psi}$ is orthogonal to $u_t$ (see \Cref{le3.1}). That is $a_0=0$. Thus,
\begin{equation}\label{e.n4}
  \hat{\psi} (x',x_N)=\sum^{\infty}_{k=1} a_k g_k(x_N) \xi_k(x').
\end{equation}
Since $\int_{\omega}\xi_k=0$ for every $k$, we get
\begin{equation*}
  \int_{\omega} \hat{\psi}_N=0.
\end{equation*}
That is, we obtain \eqref{e.new3}.
\end{proof}

Then the Fr\'{e}chet derivative $H_t$ is obtained by the following:
\begin{proposition} \label{Pr31}
For any $t\in I_{\delta}$ and any $h \in \widetilde{X}_2$, $H_{t}: \widetilde{X}_2\rightarrow \widetilde{X}_1$ is the linear continuous operator defined by
\begin{equation}\label{eq25}
H_{t}(h)(x') =t \partial_{\nu} \hat{\psi} (x',t),
\end{equation}
where $\hat{\psi}$ is the solution of \eqref{psi_h2}.
\end{proposition}

\begin{proof}
Since $F\in C^1$, we calculate the linear operator obtained by the directional derivative of $F$ with respect to $h \in \widetilde{X}_2$ at $(t,0)$. Recall that $h \in \widetilde{X}_2$ means
\begin{equation}\label{condition_w}
  h \in C^{2,\alpha}(\bar \omega), ~~\partial_{\eta} h=0~~\mbox{ and }~~\int_{\omega} h=0.
\end{equation}

Fix $t\in I_{\delta}$, let $v=t+sh$ and $(u_v,\lambda_v)$ be the unique solution of \eqref{dir}. That is, $\lambda_v$ is the first eigenvalue and $u_v$ is the positive eigenfunction with $\|u_v/\sqrt{v}\|_{L^2(\Omega_v)}=\sqrt{|\omega|/2}$ (see \Cref{dirichlet}). Recall the diffeomorphism $Y_v$ defined in \eqref{Y}. Then the function $\phi= Y_v^* u_v$ solves the problem \eqref{formula-new} with $\lambda=\lambda_v$.

\medskip

Recall now the definition of $u_t(x_N)$ given in \eqref{so:U_t} and define $\Phi= Y_v^* u_t$, which has the explicit expression:
\begin{equation} \label{bordo}
\Phi(y',y_N)=\cos \left(\sqrt{\lambda_t}v(y')y_N \right)
=\cos \left(\sqrt{\lambda_t}(t+sh(y'))y_N \right)
\end{equation}
and $\Phi$ satisfies
\begin{equation*}
     \begin{cases}
     \Delta_{g} \Phi +\lambda_t \Phi=0 &\mbox{in $\Omega_1$},\\
      \Phi=\cos \left(\sqrt{\lambda_t}(t+sh) \right) &\mbox{on $\Gamma_1$},\\
     \partial_{\eta} {\Phi} =0 &\mbox{on $\tilde \Gamma_1$}.
      \end{cases}
\end{equation*}

\medskip

Let $\psi= \phi- \Phi$ and $\mu=\lambda_v-\lambda_t$. Then $\psi$ satisfies
\begin{equation}\label{eq37}
     \begin{cases}
     \Delta_g \psi+\lambda \psi+\mu \Phi=0  &\mbox{in $\Omega_1$},\\
      \psi=-\cos \left(\sqrt{\lambda_t}(t+sh) \right) &\mbox{on $\Gamma_1$},\\
      \partial_{\eta}  \psi=0  &\mbox{on $\widetilde \Gamma_1$}.\\
      \end{cases}
    \end{equation}
 Obviously, $\psi$ is differentiable with respect to $s$. Note that $\psi\equiv 0$ if $s=0$ (c.f. \Cref{re2.3}). Let us set
\begin{equation} \label{Phi0}
\dot{\psi}=\partial_{s}\psi|_{s=0},\quad
\Phi_0(y',y_N)=\Phi(y',y_N)|_{s=0}
=\cos \left(\sqrt{\lambda_t}ty_N \right).
\end{equation}
Differentiating \eqref{eq37} with respect to $s$ at  $s=0$ and taking into account \eqref{bordo} and $\sqrt{\lambda_t}=\pi/2t$, we have
\begin{equation*}
  \begin{cases}
  \sum_{i=1}^{N-1} \dot{\psi}_{ii}+\frac{1}{t^2} \dot{\psi}_{NN}+\lambda_t \dot{\psi}+\dot{\mu} \Phi_0=\, 0  ~~&\mbox{ in }~~\Omega_1,\\
      \dot{\psi}= \, \sqrt{\lambda_t} \sin\left(\sqrt{\lambda_t}t\right) {h}=\sqrt{\lambda_t}{h}~~&\mbox{ on }~~\Gamma_1,\\
     \partial_{\eta} \dot{\psi}=\, 0~~&\mbox{ on }~~\widetilde \Gamma_1.\\
  \end{cases}
\end{equation*}
Let $(z',z_N)=(y',ty_N)$, $\hat{\psi}(z',z_N)=\dot{\psi}(y',y_N)/t$ and $\hat{\Phi}(z',z_N)=\Phi_0(y',y_N)/t$. Then $\hat{\psi}$ satisfies
\begin{equation}\label{psi_h}
  \begin{cases}
     \Delta \hat{\psi}+\lambda_t \hat{\psi}+\dot{\mu} \hat{\Phi}=\, 0  ~~&\mbox{ in }~~\Omega_t,\\
      \hat{\psi}=\frac{1}{t}\sqrt{\lambda_t} {h}= \frac{\pi}{2t^2} h~~&\mbox{ on }~~\Gamma_t,\\
     \partial_{\eta} \hat{\psi}=\, 0~~&\mbox{ on }~~\widetilde \Gamma_t.\\
  \end{cases}
\end{equation}
From \eqref{Phi0}, we obtain
\begin{equation}\label{hat_Phi_ex}
  \hat{\Phi}(z',z_N)=\frac{\Phi_0(y',y_N)}{t}=\frac{\cos \left(\sqrt{\lambda_t}ty_N \right)}{t}
  =\frac{\cos \left(\sqrt{\lambda_t}z_N \right)}{t}
\end{equation}
and have that $\hat{\Phi}$ satisfies
\begin{equation}\label{hat_Phi}
     \begin{cases}
     \Delta \hat \Phi +\lambda_t \hat \Phi=0 &\mbox{in $\Omega_t$},\\
     \hat \Phi=0 &\mbox{on $\Gamma_t$},\\
     \partial_{\eta} {\hat \Phi} =0 &\mbox{on $\widetilde \Gamma_t$}.
      \end{cases}
\end{equation}
Next, we multiply the first equation of \eqref{psi_h} by $\hat{\Phi}$ and integrate it over $\Omega_t$. Then, by considering the boundary condition of \eqref{psi_h} and \eqref{hat_Phi} and the equation in \eqref{hat_Phi}, we obtain
\begin{equation*}
  \begin{aligned}
    0&=\int_{\Omega_t} \Delta \hat{\psi} \cdot \hat{\Phi}+\lambda_t \hat{\psi} \hat{\Phi}
    +\dot{\mu} \hat{\Phi}^2\\
    &=-\int_{\Omega_t} \nabla \hat{\psi}\cdot \nabla \hat{\Phi}+\int_{\Gamma_t}\langle \nabla \hat{\psi},\nu \rangle \hat{\Phi}+\int_{\widetilde \Gamma_t}\langle \nabla \hat{\psi},\eta \rangle \hat{\Phi}
    +\int_{\Omega_t} \lambda_t \hat{\psi} \hat{\Phi}
    +\dot{\mu} \hat{\Phi}^2\\
    &=-\int_{\Omega_t} \nabla \hat{\psi}\cdot \nabla \hat{\Phi}
    +\int_{\Omega_t} \lambda_t \hat{\psi} \hat{\Phi}
    +\dot{\mu} \hat{\Phi}^2\\
    &=\int_{\Omega_t} \hat{\psi}\cdot \Delta \hat{\Phi}-\int_{\Gamma_t}\langle \nabla \hat{\Phi},\nu \rangle \hat{\psi}-\int_{\widetilde \Gamma_t}\langle \nabla \hat{\Phi},\eta \rangle \hat{\psi}
    +\int_{\Omega_t} \lambda_t \hat{\psi} \hat{\Phi}
    +\dot{\mu} \hat{\Phi}^2\\
    &=\int_{\Omega_t} \hat{\psi}\cdot \left(\Delta \hat{\Phi}+\lambda_t \hat{\Phi}\right)
    -\int_{\Gamma_t}\langle \nabla \hat{\Phi},\nu \rangle \hat{\psi}
    +\int_{\Omega_t} \dot{\mu} \hat{\Phi}^2\\
    &=-\int_{\Gamma_t}\langle \nabla \hat{\Phi},\nu \rangle \hat{\psi}
    +\int_{\Omega_t} \dot{\mu} \hat{\Phi}^2.
  \end{aligned}
\end{equation*}
By considering the fact $\int_{\omega} h=0$ (see \eqref{condition_w}), the boundary condition of \eqref{bordo} and \eqref{hat_Phi_ex}, we have
\begin{equation*}
  \begin{aligned}
    \int_{\Gamma_t}\langle \nabla \hat{\Phi},\nu \rangle \hat{\psi}
    =\int_{\Gamma_t}\hat{\Phi}_N \hat{\psi}
    =-\int_{\Gamma_t}\frac{\sqrt{\lambda_t}}{t}\sin\left(\sqrt{\lambda_t}t\right)
    \cdot \left(\frac{\sqrt{\lambda_t}}{t} h\right)
    =-\frac{\pi}{2t^3} \int_{\omega} h
    =0.
  \end{aligned}
\end{equation*}
Then
\begin{equation*}
  \int_{\Omega_t} \dot{\mu} \hat{\Phi}^2=0.
\end{equation*}
Thus, we conclude that $\dot{\mu}=0$ and \eqref{psi_h} becomes \eqref{psi_h2}.
In other words,
\begin{equation} \label{sanremo}
     \dot{\psi}(y',y_N)= t \hat{\psi}(z',z_N)= t \hat{\psi}(y',ty_N).
\end{equation}

\medskip

In addition, set
\begin{equation*}
  \dot{\phi}=\partial_{s}\phi|_{s=0}~~\mbox{ and }~~\dot{\Phi}=\partial_{s}\Phi|_{s=0}.
\end{equation*}
Then by noting \eqref{bordo} we have
\begin{equation*}
  \begin{aligned}
    \dot{\phi}(y',y_N)=&\dot{\psi}(y',y_N)+\dot{\Phi}(y',y_N)\\
    =&t \hat{\psi}(y',ty_N)+\dot{\Phi}(y',y_N)\\
    =& t \hat{\psi}(y',ty_N)
         -\sin \left(\sqrt{\lambda_t}t y_N \right)\sqrt{\lambda_t} y_N h\\
    =&t \hat{\psi}(y',ty_N)
         -\frac{\pi}{2t} y_N\sin \left(\frac{\pi}{2} y_N \right) h.
  \end{aligned}
\end{equation*}
Hence, the derivative of $\dot{\phi}$ with respect to $y_N$ is
\begin{equation}\label{phi_N}
  \dot{\phi}_N(y',y_N)
  =t^2 \hat{\psi}_N(y',ty_N)
         -\frac{\pi}{2t} \sin \left(\frac{\pi}{2} y_N \right) h
         -\frac{\pi^2}{4t} y_N \cos \left(\frac{\pi}{2} y_N \right) h.
\end{equation}

\medskip

On the other hand, we can write $\phi$ as a function of $s$ as
\begin{equation*}
  \begin{aligned}
    \phi(s)=&\phi(0) +s \dot{\phi} +\mathcal{O}(s^{2})\\
    =&\Phi(0) +s \dot{\phi} +\mathcal{O}(s^{2})\\
    =&\cos \left(\sqrt{\lambda_t}(t+sh)y_N \right) |_{s=0}+s \dot{\phi} +\mathcal{O}(s^{2})\\
    =&\cos \left(\sqrt{\lambda_t}ty_N \right)+s \dot{\phi} +\mathcal{O}(s^{2}).
  \end{aligned}
\end{equation*}
Take derivatives of $\phi$ with respect to $y_i$ and $y_N$,
\begin{equation}\label{phi_i_N}
  \begin{aligned}
    \phi_i(s)=&s\dot{\phi}_i+\mathcal{O}(s^{2}),\\
    \phi_N(s)=&-\sqrt{\lambda_t}t \sin \left(\sqrt{\lambda_t
    }ty_N \right)+s \dot{\phi}_N +\mathcal{O}(s^{2})\\
    =&-\frac{\pi}{2} \sin \left(\frac{\pi}{2} y_N \right)+s \dot{\phi}_N +\mathcal{O}(s^{2}).
  \end{aligned}
\end{equation}
If $s=0$, we obtain $\phi_i=0$ and $\phi_N=-\frac{\pi}{2} \sin \left(\frac{\pi}{2} y_N \right)$.

\medskip

Note that
\begin{equation*}
  |\nabla u_v|^2=g^{ij}\phi_i\phi_j.
\end{equation*}
Then we have
\begin{equation}\label{u-nu-1}
  \begin{aligned}
    D_s (\partial_{\nu} u_v)
    &:=\left(\frac{\partial u_v}{\partial \nu}\right)_s \bigg|_{s=0,x_N=t}\\
    &=-\left(\sqrt{g^{ij}\phi_i\phi_j}\right)_s \bigg|_{s=0,y_N=1}\\
    &=-\frac{1}{2}\left(g^{ij}\phi_i\phi_j\right)^{-\frac{1}{2}}\cdot\left(g^{ij}\phi_i\phi_j\right)_s \big|_{s=0,y_N=1}.
  \end{aligned}
\end{equation}
Observe that if $i\neq N$, by \eqref{phi_i_N}, we have $\phi_i(s)=\mathcal{O}(s)$. Thus, $g^{ij}\phi_i\phi_j=\mathcal{O}(s^2)$ for $i,j\neq N$ and then we have
\begin{equation*}
  \left(g^{ij}\phi_i\phi_j\right)_s \big|_{s=0}=0.
\end{equation*}
Additionally, by \eqref{g-inverse}, we have for $i\neq N$,
\begin{equation*}
  g^{iN}=-\frac{v_iy_N}{v}=-\frac{s h_i y_N}{t+sh}=\mathcal{O}(s).
\end{equation*}
Similarly, we obtain $g^{iN}\phi_i\phi_N=\mathcal{O}(s^2)$ and
\begin{equation*}
  \left(g^{iN}\phi_i\phi_N\right)_s \big|_{s=0}=0.
\end{equation*}
Therefore,
\begin{equation}\label{u-nu-12}
    D_s (\partial_{\nu} u_v)
    =-\frac{1}{2}\left(g^{NN}\phi_N^2\right)^{-\frac{1}{2}}\cdot\left(g^{NN}\phi_N^2\right)_s \big|_{s=0,y_N=1}.
\end{equation}
Recall \eqref{g-inverse} again,
\begin{equation*}
  g^{NN}=\frac{1+|Dv|^2y_N^2}{v^2}=\frac{1+s^2|Dh|^2y_N^2}{(t+sh)^2}.
\end{equation*}
Then \eqref{u-nu-12} becomes
\begin{equation*}
  \begin{aligned}
    D_s (\partial_{\nu} u_v)
    &=-\frac{1}{2}\frac{t}{|\phi_N|}\left(\frac{2\phi_N \dot \phi_N}{t^2}-\phi_N^2\frac{2h}{t^3}\right) \bigg|_{s=0}\\
    &=\frac{\dot \phi_N}{t}+\frac{\pi}{2t^2}h
  \end{aligned}
\end{equation*}
on $\Gamma_1$.
By combining \eqref{phi_N}, we have
\begin{equation*}
  D_s (\partial_{\nu} u_v)=t \hat{\psi}_N(y',t).
\end{equation*}
Next, by the definition of $F$ (see \eqref{F}) and $H_t$ (see \eqref{Ht}), we know
\begin{equation*}
\begin{aligned}
  H_t(h)(y')
  &=D_s (\partial_{\nu} u_v)-\frac{1}{|\omega|}\int_{\omega} D_s (\partial_{\nu} u_v)\\
  &=t \hat{\psi}_N(y',t)-\frac{t}{|\omega|}\int_{\omega}\hat{\psi}_N(y',t).
  \end{aligned}
\end{equation*}
Then by \Cref{le.new}, we get \eqref{eq25}.

\end{proof}

\section{Study of the linearized operator $H_t$}
\label{sec:linearized_operator_H_t}

In this section, we study the properties of the operator $H_t$ given by \eqref{eq25}. As in the proof of \Cref{le.new}, let $\{\xi_k\}_{k\geq 0}$ be the sequence of normalized eigenfunctions (i.e., $\|\xi_k\|_{L^2(\omega)}=1$) of the Laplacian operator on $\omega \subset \mathbb{R}^{N-1}$ with Neumann boundary conditions and $\sigma_k$ ($k\geq1$) are the corresponding sequence of eigenvalues.
We will show that:
\begin{proposition} \label{H}
For a fixed $t_*$ and $t\in I_{\delta}=(t_*-\delta,t_*+\delta)$, the operator $H_t$ is an essentially self-adjoint Fredholm operator of index $0$ and has a sequence of eigenvalues denoted by $\left\{\mu_{t,k}\right\}_{k\geq 1}$. Let
\begin{equation*}
  t_k=\frac{\pi}{2\sqrt{\sigma_k}},\quad k\geq 1.
\end{equation*}
Then the eigenvalues $\left\{\mu_{t,k}\right\}_{k\geq 1}$ have the following expressions:
\begin{equation} \label{eigenvalues}
\mu_{t,k}=
\begin{cases}
  -\frac{\pi^2}{4t^2}\sqrt{1-\frac{t^2}{t_k^2}}
  \tan \left(\frac{\pi}{2}\sqrt{1-\frac{t^2}{t_k^2}} \right),~~~~&t<t_k,\\
  0,~~~~&t=t_k,\\
  \frac{\pi^2}{4t^2} \sqrt{\frac{t^2}{t_k^2}-1}
  \left(1-\frac{2}{e^{\pi \sqrt{t^2/t_k^2-1} }+1}\right),~~~~&t>t_k.
  \end{cases}
\end{equation}
Moreover, the eigenfunction associated to $\mu_{t,k}$ is just $\xi_k$ for any $k\geq 1$.
\end{proposition}

\begin{remark}\label{propo4.1}
We remark here that later $t_*$ will be chosen as $t_j=\frac{\pi}{2\sqrt{\sigma_j}}$ (see \Cref{sec:proof_Tmain}) where $\sigma_j$ is a simple eigenvalue, therefore \eqref{eigenvalues} could just be written as $t_j<t_k, t_j=t_k$ and $t_j>t_k$ if $\delta$ is sufficiently small.
\end{remark}

\begin{proof}
Recall that $\hat \psi$ satisfies \eqref{psi_h2} and $w$ can be written as the Fourier series shown in \eqref{e.new}. Then we have \eqref{e.n4}. Inserting \eqref{e.n4} into \eqref{psi_h2}, we have \eqref{psi_h3}.

\medskip

Next, the solution of \eqref{psi_h3} can be discussed as following:\\
\textbf{(i).} If $\lambda_t=\left(\frac{\pi}{2t}\right)^2=\sigma_k$, i.e., $t=t_k$, we have $g_k''=0$. Thus,
\begin{equation}\label{g_k1}
g_k\equiv \frac{\pi}{2t^2}
~~\mbox{ and }~~
g_k'\equiv 0.
\end{equation}
\textbf{(ii).} If $\lambda_t>\sigma_k$, i.e., $t<t_k$, we have
\begin{equation*}\label{g_k2_1}
g_k(x_N)=\frac{\pi}{2t^2} \frac{\cos\left(\sqrt{\lambda_t-\sigma_k} x_N\right)}{\cos\left(\sqrt{\lambda_t-\sigma_k} t\right)}
\end{equation*}
and
\begin{equation}\label{g_k2_2}
g_k'(x_N)=-\frac{\pi}{2t^2}
\frac{\sqrt{\lambda_t-\sigma_k} \sin\left(\sqrt{\lambda_t-\sigma_k} x_N\right)}{\cos\left(\sqrt{\lambda_t-\sigma_k} t\right)}.
\end{equation}
\textbf{(iii).} If $ \lambda_t<\sigma_k$, i.e., $t>t_k$, we obtain
\begin{equation*}\label{g_k3_1}
g_k(x_N)=\frac{\pi}{2t^2}
\frac{e^{\sqrt{\sigma_k-\lambda_t} x_N}
+e^{-\sqrt{\sigma_k-\lambda_t} x_N}}{e^{\sqrt{\sigma_k-\lambda_t} t}
+e^{-\sqrt{\sigma_k-\lambda_t} t}}
\end{equation*}
and
\begin{equation}\label{g_k3_2}
g_k'(x_N)=\frac{\pi}{2t^2} \frac{\sqrt{\sigma_k-\lambda_t}
\left(e^{\sqrt{\sigma_k-\lambda_t} x_N}
-e^{-\sqrt{\sigma_k-\lambda_t} x_N}\right)}{e^{\sqrt{\sigma_k-\lambda_t} t}
+e^{-\sqrt{\sigma_k-\lambda_t} t}}.
\end{equation}
Therefore, for any $k\geq 1$,
\begin{equation}\label{g_1pri}
  g_k'(t)=
  \begin{cases}
  -\frac{\pi^2}{4t^3}\sqrt{1-\frac{t^2}{t_k^2}}
  \tan \left(\frac{\pi}{2}\sqrt{1-\frac{t^2}{t_k^2}} \right),~~~~&t<t_k,\\
  0,~~~~&t=t_k,\\
  \frac{\pi^2}{4t^3} \sqrt{\frac{t^2}{t_k^2}-1}
  \left(1-\frac{2}{e^{\pi \sqrt{t^2/t_k^2-1} }+1}\right),~~~~&t>t_k.
  \end{cases}
\end{equation}

\medskip

Note that
\begin{equation}\label{hat_psi-N}
H_t(w)=t\partial_{\nu} \hat{\psi} (x',t)=t\hat{\psi}_N (x',t)=t\sum^{\infty}_{k=1} a_k g_k'(t) \xi_k(x').
\end{equation}
Hence, we have
\begin{equation*}
  \mu_{t,k}=tg_k'(t).
\end{equation*}
By noting \eqref{g_1pri}, we have the expressions \eqref{eigenvalues} for $\mu_{t,k}$.

Note that $t_k \to 0$ if $k \to +\infty$. So for any fixed $t$, we have $t>t_k$ for $k$ large enough and then by \eqref{g_1pri},
\begin{equation}\label{e4.1}
|g'_k(t)|\leq \frac{C}{t_k}\leq C\sqrt{\sigma_k}
\end{equation}
for some constant $C$ independent of $k$. If $w\in \widetilde{H}^1(\omega)$, by noting \eqref{e.new}, we have
\begin{equation*}
  \|\nabla w\|^2_{L^2(\omega)}=\sum^{\infty}_{k=1} a_k^2 \|\nabla \xi_k\|^2_{L^2(\omega)}
  =\sum^{\infty}_{k=1} a_k^2 \sigma_k.
\end{equation*}
By combining with \eqref{hat_psi-N} and \eqref{e4.1},
\begin{equation*}
  \|H_t(w)\|^2_{L^2(\omega)}=t^2\sum^{\infty}_{k=1} a_k^2 (g'_k(t))^2
  \leq t^2\sum^{\infty}_{k=1} a_k^2 \sigma_k
  =t^2\|\nabla w\|^2_{L^2(\omega)}.
\end{equation*}
Thus, we have $H_t(w)\in L^2(\omega)$ with $\int_{\omega} w=0$ (similarly, we denote $ H_t(w)\in \widetilde{L}^2(\omega)$ for simplicity). In conclusion, $H_t$ can be extended to the Sobolev framework:
\begin{equation*}
  H_t: \widetilde{H}^1(\omega) \to \widetilde{L}^2(\omega),
\end{equation*}
where $H_t$ is a self-adjoint operator of order $1$. By H\"{o}lder regularity, $H_t$ is a Fredholm operator of index 0.
\end{proof}

\section{Bifurcation argument and proof of the main result}
\label{sec:proof_Tmain}

In this section we conclude the proof of our main result by means of the classical Crandall-Rabinowitz Theorem. For convenience of the reader we recall its statement below, in a $C^2$ version.

\begin{theorem} \label{CR} \textbf{\mbox{(Crandall-Rabinowitz Bifurcation Theorem)}}
	Let $X$ and $Y$ be Banach spaces, and let $U\subset X$ and $I\subset\mathbb{R}$ be open domains, where we assume $0\in U$. Denote the elements of $U$ by $v$ and the elements of $I$ by $t$. Let $G:I \times U \rightarrow Y$ be a $C^{2}$ operator such that
	\begin{itemize}
		\item[i)] $G(t,0)=0$ for all $t\in I,$
		\item[ii)] $\emph{Ker~} D_{v}G(t_{*},0)=\mathbb{R}\,w$ for some $t_{*}\in I$ and some $w\in X\setminus\{0\};$
		\item[iii)] $\emph{codim Im~} D_{v}G(t_{*},0)=1;$
		\item[iv)] $ D_{t}D_{v}G(t_{*},0)(w)\notin \emph{Im~} D_{v}G(t_{*},0).$
	\end{itemize}
	Then there exists a nontrivial $C^1$ curve
	\begin{equation}\label{eq501}
		(-\e, \e) \ni s \mapsto (t(s), v(s)) \in I \times X,
	\end{equation}
for some $\e>0$, such that:
\begin{enumerate}
	\item $t(0)=t_*$, $t'(0)=0$, $v(0)=0$, $v'(0)=w$.
	\item $G\big(t(s),v(s)\big)=0$ for all $ s\in(-\delta,+\delta)$.
\end{enumerate}
	Moreover, there exists a neighborhood $\mathcal{N}$ of $(t_{*},0)$ in $X \times \Gamma$ such that all solutions of the equation $G(t,v)=0$ in $\mathcal{N}$ belong to the trivial solution line $\{(t,0)\}$ or to the curve (\ref{eq501}). The intersection $(t_{*},0)$ is called a bifurcation point.
\end{theorem}

We are now in a position to give the~\\
\noindent\textbf{Proof of \Cref{Tmain}.} For this we only need to apply \Cref{CR} to our operator $F$. If $\sigma_j$ is simple, let
\begin{equation*}
t_{*}=t_j=\frac{\pi}{2\sqrt{\sigma_j}}.
\end{equation*}
Next, set
\begin{equation*}
I=(t_{*}/2,2t_{*}), \quad X=\widetilde{X}_2, \quad U=\left\{w \in\widetilde{X}_2: \|w\|_{L^{\infty}(\omega)}<t_{*}/2 \right\}, \quad Y=\widetilde{X}_1.
\end{equation*}

Note that $I\times U\subset X_2^+$ and hence by \Cref{dirichlet}, $F:I\times U\to Y$ is a smooth map. Next, since
\begin{equation*}
F(t,0)=\frac{\partial u_t}{\partial {\nu}} (x', t) - \frac{1}{|\omega|} \, \int_{\omega} \, \frac{\partial u_t}{\partial \nu}=0,~\forall ~t\in I,
\end{equation*}
i) is satisfied.

Since $\sigma_j$ is simple, according to \eqref{g_1pri} and \eqref{hat_psi-N}, the kernel of $H_{t_*}$ is spanned by $\xi_j$. So ii) is satisfied with $w= \xi_j$.

Since $H_t$ is essentially self-adjoint, we have that:
$$ \emph{Im~}D_w F(t_*,0) = \{z \in \widetilde X_1:\ \int_{\omega} z \xi_j =0\},$$
which has codimension 1, as required by iii).

It remains to check condition iv), i.e. the so-called {\it transversality condition}. Recall the linearization of the normal derivative operator:
\begin{equation*}
  D_w F(t,w)|_{w=0}=H_t (w)=t\sum^{\infty}_{k=1} a_k g_k'(t) \xi_k(x'),
\end{equation*}
where $g'_k(t)$ is shown in \eqref{g_1pri}. Next, we directly compute the derivatives at $t=t_*$ to obtain:
\begin{equation*}
D_t D_{w}F(t_*,0)(\xi_j)=\frac{\pi^3}{4 t_*^3}\xi_j \in \mathrm{Ker}~H_{t_*}.
\end{equation*}

Thus, we conclude that $D_t D_{w}F(t_*,0)(\xi_j)$ does not belong to $ \emph{Im~} D_w F(t_*,0)$. Therefore, all assumptions of \Cref{CR} are satisfied. By \Cref{CR}, the conclusion of \Cref{Tmain} follows. ~\qed~\\

\medskip

\printbibliography

\end{document}